\documentclass[11pt,reqno]{amsart}
\usepackage{latexsym,amssymb}

\usepackage{cite}

\newtheorem{defi}{Definition}[section]
\newtheorem{proposition}[defi]{Proposition}

\newtheorem{rem}[defi]{Remark}
\newtheorem{lemma}[defi]{Lemma}
\newtheorem{theorem}[defi]{Theorem}

\newcommand{\defeq}{\mathrel{\mathrm{\raise0.1ex\hbox{:}\hbox{=}\strut}}}

\def\dd{\:\mathrm{d}}

\def\R{\mathbb R}
\def\NN{\mathbb{N}}
\def\E{\mathbb E}
\def\EE{\mathbb E}
\def\ee{\mathrm{e}}

\def\ve{\varepsilon}

\def\la{\langle}
\def\ra{\rangle}

\def\ein{\mathbf{1}}

\begin{document}

\numberwithin{equation}{section}

\title[Invariant measures for monotone SPDE's ]{Invariant measures for monotone SPDE's with multiplicative noise term }
\author[A. Es--Sarhir]{Abdelhadi Es--Sarhir }
\address{D\'epartement de Math\'ematiques, Facult\'e des Sciences, Universit\'e Ibn Zohr\\ 
BP. 8106, Cit\'e Dakhla, 80000 Agadir, Morocco}
\email{a.es-sarhir@uiz.ac.ma}
\author[M. Scheutzow]{Michael Scheutzow}
\author[J. M. T\"olle]{Jonas M. T\"olle}
\address{Technische Universit\"at Berlin, Fakult\"at II, Institut f\"ur Mathematik, Sekr. MA 7-5\newline \noindent Stra{\ss}e des 17. Juni 136, D-10623 Berlin, Germany}
\email{ms@math.tu-berlin.de}
\email{jonasmtoelle@gmail.com}
\author[O. van Gaans]{Onno van Gaans}
\address{Universiteit Leiden, Mathematisch Instituut\\
 Postbus 9512, 2300 RA Leiden, The Netherlands}
\email{vangaans@math.leidenuniv.nl}


\begin{abstract} We study
diffusion processes corresponding to infinite dimensional semilinear
stochastic differential equations with local Lipschitz drift term
and an arbitrary Lipschitz diffusion coefficient. We prove tightness
and the Feller property of the solution to show existence of an
invariant measure.
As an application we discuss stochastic reaction
diffusion equations.
\end{abstract}

\maketitle
\section{Introduction and preliminaries}
\noindent We are dealing with the following semilinear stochastic
differential equation on a real separable Hilbert space $H$
\begin{equation}\label{sde0}
 \left\{
\begin{array}{ll}
du(t)=\Big(Au(t)+F(u(t))\Big)dt+B(u(t))dW_t,\quad\mbox{for}\:\: t\geq0,\\
u(0)=x\in H,
\end{array}
\right.
\end{equation}
 where $A$ is a self adjoint operator with negative type
$\omega$ on $H$ and compact resolvent $A^{-1}$, $F\colon\,H\to H$ is
a continuous nonlinear mapping, and $(W_t)_{t\ge 0}$ is a cylindrical Wiener
process in a separable real Hilbert space $U$ defined on a filtered
probability space $(\Omega,\mathcal{F},(\mathcal{F}_t)_{t\ge
0},\mathbb{P})$. The coefficient $B$ maps $H$ into the space of
Hilbert-Schmidt operators $\mathcal{L}_{HS}(U,H)$ from $U$ into $H$
and is assumed to be Lipschitz from $H$ into $\mathcal{L}_{HS}(U,H)$
with Lipschitz constant $L$.

\noindent Equation \eqref{sde0} can been seen as an abstract
formulation of reaction diffusion equations perturbed by random
noise. In this model of equation the nonlinear drift $F$ is locally
Lipschitz and has additional dissipative properties. This special
structure of $F$ has been used to analyze \eqref{sde0} in a space of
continuous functions as a subspace of $L^2$ (see \cite{cerrai:03}).
It is our main aim to analyze equation \eqref{sde0} for locally
Lipschitz $F$ with suitable quasi-dissipative properties in a
general Hilbert space setting. We are mainly interested in the
existence of an invariant measure for \eqref{sde0} without condition
on the Lipschitz constant of $B$. Our analysis is based on a
Lyapunov type assumption on the coefficient $F$ and the compactness
of the linear part. For a general theory of reaction diffusion
equations including both cases of additive and multiplicative noise
perturbations we refer to the monographs \cite{DaZa:96,cerrai00} and
the works \cite{cerrai99, cerrai:03, cerrai05,GaGo:94,
GaGo:97,GoMa,Ma-Za}. We note that the global mild solution $u$ satisfies
the following integral equation
\begin{equation}\label{S}
u(t)=e^{tA}x+ \int_0^t e^{(t-s)A}F(u(s))ds+\int_0^t e^{(t-s)A}
B(u(s))dW_s,\quad t\ge 0,
\end{equation}
with transition semigroup
$$
P_t\varphi(x)=\E(\varphi(u(t,x))),\quad x\in H,\:\:t\geq0,
$$
defined on the space of all bounded measurable functions on $H$. An
\emph{invariant measure} for (\ref{sde0}) is a Borel probability
measure $\mu$ on $H$ such that
\begin{equation*}
P_t^*\mu=\mu\quad\mbox{for all }t\ge 0,
\end{equation*}
where $P_t^*$ denotes the adjoint of $P_t$.

\noindent In the literature there are several conditions ensuring
the existence of such measures $\mu$, one of them is based on
Krylov-Bogoliubov's theorem using a compactness property
for the underlying semigroup generated by the linear part in
\eqref{sde0} and  boundedness in probability of solutions. However,
to have the latter property is not in general straightforward: In
most cases one checks the boundedness of the moments of solutions
which requires some specific conditions on the coefficients $A$, $F$
and $B$. In \cite{DGZ} it was proved that if the coefficients of
\eqref{sde0} satisfy a dissipativity condition then \eqref{sde0} has
a bounded solution which has an invariant measure by using so called
\emph{remote start method}. However this dissipativity assumption is
strong in the sense that the Lipschitz constants of $F$ and $B$
should be small compared to the exponential growth of the semigroup
generated by $A$. Of course in the use of the compactness argument
the dissipativity on the term $B$ can be relaxed and one can suppose
the boundedness of $B$ or the existence of a bounded solution to
show the existence of an invariant measure. Our hypothesis ${\bf(
H_4)}$ on the drift $F$ (see below) is inspired by \cite{Es-Ga-Sch},
which discusses existence of an invariant measure for stochastic
delay equations in finite dimensions. It turns out that our
condition on $F$ allows general terms $B$ which are only Lipschitz.
Let us now define the following interpolation spaces.

\noindent For $\gamma\in\R$ let
$$V_\gamma\defeq (D((-A)^{\gamma}),\|\cdot\|_{\gamma}),\quad\mbox{
where $\la x,y\ra_{\gamma} = \la (-A)^{\gamma} x,(-A)^{\gamma} y\ra$
for $x,\: y\in V_{\gamma}$}.$$ Note that, since $A$ has a compact
resolvent, the embedding $V_{\gamma}\hookrightarrow H$ is compact
for $\gamma>0$. In the following $\|\cdot\|_0$ and  $\|\cdot\|_{HS}$
denote the $H$-norm and the Hilbert-Schmidt operator norm
respectively. We shall formulate our assumptions:

\begin{enumerate}

\item [${\bf (H_0)}$]

$A$ is selfadjoint and $\|e^{tA}\|\leq e^{-\omega t}$ for some
$\omega>0$.
\item [${\bf (H_1)}$] $F:\:(E,\|\cdot\|_E)\longrightarrow E$ is locally Lipschitz continuous
and bounded on bounded sets of the Banach space $E\subset H$. The
part of $A$ in $E$ denoted by $A_{E}$ generates a strongly
continuous semigroup on $E$ and $E$ is a Banach space continuously,
densely, and as a Borel subset embedded in $H$. The embedding
$V_{\gamma}\hookrightarrow E$ is continuous, $\gamma\in (\frac
14,\frac 12)$.
\item[${\bf (H_2)}$] There exists an increasing function
$a:\:\R^+\longrightarrow\R^+$ such that
\begin{equation}\label{a}
_{E}\la F(y+z),y^{\ast}\ra_{E^{\ast}} \leq
a(\|z\|_E)(1+\|y\|_E)\quad\mbox{for}\:\:y,\:z\in E,\:y^{\ast}\in
\partial\:\|y\|_E,
\end{equation}
where $\partial\|y\|_E$ denotes the subdifferential of
$\left\|\cdot\right\|_E$ at $y$.
\item[${\bf (H_3)}$] There exists $\kappa>0$ such that
$$
\la F(u)-F(v),u-v\ra\leq \kappa \|u-v\|_0^2,\quad u,\:v\in E.
$$

\end{enumerate}

\noindent We end this introduction by the following definition.

\begin{defi}\label{Definition}
A mild solution of equation \eqref{sde0} is an
$\mathcal{F}_t$-adapted process $u$ such that $u\in
C([0,+\infty),E)$ a.s.\ and which satisfies the following integral
equation
\begin{equation}\label{mild-solution}
u(t)=e^{tA}x + \int_0^t e^{(t-s)A}F(u(s))ds+\int_0^t
e^{(t-s)A}B(u(s))dW_s,\quad t\ge 0.
\end{equation}
\end{defi}
\section{Existence and uniqueness of solutions}

\noindent In this section we show existence of a unique global
solution $(u(t))_{t\geq 0}$ for the equation \eqref{sde0}. We start
by the following lemma. For a proof we refer to \cite{DZ, Ga} (see
also the proof of Theorem 2.3 in \cite{Seidler}.)
\begin{lemma}\label{lemma} Let $p>2$ and $\eta\colon [0,T]\times\Omega\to \mathcal{L}_{HS}(U,H)$
 be a progressively measurable process with
\[
\EE\int_0^T \|\eta(s)\|^p_{HS}\,\mathrm{d}s<\infty.
\]
If $\gamma+\frac 1p<\frac12$, then  $\int_0^t
e^{(t-s)A}\eta(s)\,\mathrm{d} W(s)$ has a continuous version in
$V_\gamma$.
% and there exists a constant $\kappa_p(T,\gamma)$ with
%$\kappa_p(T,\gamma)\xrightarrow[T\to 0]{}0$ such that
%\begin{equation}\label{stochconvestim} \EE\sup_{0\le t\le T}
%\left\|\int_0^t e^{(t-s)A} \eta(s)
%\,\mathrm{d}W_s\right\|_{\gamma}^p \le \kappa_p(T,\gamma) \EE
%\int_0^T \|\eta(s)\|_{HS}^p\,\mathrm{d}s.
%\end{equation}

\end{lemma}
%\begin{proof}
% Define
%$$ Y(t)\defeq \int_0^t(t-s)^{-\alpha}e^{(t-s)A}\eta(s)\:dW_s.
%$$
%Using the factorization formula as in \cite{DZ, Ga} we can write
%$$
%\int_0^t
%e^{(t-s)A}\eta(s)\:dW_s=\frac{\sin\pi\alpha}{\pi}R_{\alpha}Y(t)
%$$
%where
%$$
%R_{\alpha}f(t)=\int_0^t(t-s)^{\alpha-1}e^{(t-s)A}f(s)\:ds,\quad\mbox{for}\:\:
%\gamma+\frac 1p <\alpha,\:\:\gamma+\alpha<\frac12,
%$$
%defines a bounded linear operator from $L^p([0,T],H)$ into
%$C([0,T],H)$. Therefore, we can write

%\begin{equation*}
%\begin{split}
%\E\left(\sup\limits_{0\leq t\leq T}\|\int_0^t
%e^{(t-s)A}\eta(s)\:dW_s\|_{\gamma}^p\right)^{\frac
%1p}&=\E\left(\sup\limits_{0\leq t\leq
%T}\|\frac{\sin\pi\alpha}{\pi}R_{\alpha}(-A)^{\gamma}Y(t)\|_0^p\right)^{\frac 1p}\\
%&\leq
%\frac{1}{\pi}\|R_{\alpha}\|\left(\E\|(-A)^{\gamma}Y(\cdot)\|_{L^p([0,T],H)}^p\right)^{\frac
%1p}.
%\end{split}
%\end{equation*}
%Using Burkholder-Davis-Gundy's inequality we obtain
%\begin{equation*}
%\begin{split}
%\E\|(-A)^{\gamma}Y\|_{L^p([0,T],H)}^p&=\E\int_0^T\|(-A)^{\gamma}Y(t)\|_0^p\:dt\\
%&=\int_0^T\E\|\int_0^t(t-s)^{-\alpha}(-A)^{\gamma}e^{(t-s)A}\eta(s)\:dW_s\|_0^p\:dt\\
%&\leq
%c_p\E\int_0^T\left(\int_0^t(t-s)^{-2\alpha}\|(-A)^{\gamma}e^{(t-s)A}\eta(s)\|_{HS}^2\:ds\right)^{\frac
%p2}\:dt\\
%(\mbox{Young's inequality} )&\leq
%c_p\left(\int_0^Ts^{-2(\alpha+\gamma)}\:ds\right)^{\frac
%p2}\cdot \E\int_0^T\|\eta(s)\|^p_{HS}\:ds\\
%&\leq c_p(\alpha,\gamma) T^{\frac
%p2(1-2(\alpha+\gamma))}\cdot\E\int_0^T\|\eta(s)\|^p_{HS}\:ds,
%\end{split}
%\end{equation*}
%which concludes the proof of the lemma.
%\end{proof}
\begin{theorem}\label{existence}
Under hypotheses ${\bf (H_0)}$, ${\bf (H_1)}$,  ${\bf (H_2)}$ and
${\bf (H_3)}$ equation \eqref{sde0} has a unique global mild
solution for each initial condition $x\in E$.
\end{theorem}
\begin{proof}
\noindent For $T>0$, $p>4$, and an $E$-valued, progressively
measurable process $v$ with
$$\E\int_0^T\|v(s,\omega)\|_0^p\:ds<+\infty$$
we introduce on $[0,T]$ the following differential equation
\begin{equation}\label{sde1}
 \left\{
\begin{array}{ll}
dz(t)=&\Big(Az(t)+F(z(t))\Big)dt+B(v(t))dW_t,\quad
t\in [0,T],\\
z(0)\:=&x \in E.
\end{array}
\right.
\end{equation}
\noindent We remark that since we assumed in ${\bf ( H_1)}$ that the
embedding $V_{\gamma}\hookrightarrow E$ is continuous we have by
Lemma \ref{lemma} that the stochastic convolution $\int_0^t
e^{(t-s)A} B(v(s))\:dW_s$ has a continuous version in $E$. Hence by
using hypothesis ${\bf ( H_0)}$-${\bf ( H_2)}$ and Theorem 7.10 in
\cite{DPZ1}, equation \eqref{sde1} has a unique mild solution
$z$ with paths in $C([0,+\infty), E)$. We now introduce the space
$\mathcal{K}$ of all $H$-valued predictable processes $z$ defined on
the interval $[0,T]$ such that
$$
\|z\|_{\mathcal{K}}=\sup_{0\le t\le T}\E\left(
\|z(t)\|_0^p\right)^{1/p}<\infty.
$$

\noindent Clearly, $\|\cdot\|_{\mathcal{K}}$ is a norm on
$\mathcal{K}$ and $(\mathcal{K},\|\cdot\|_{\mathcal{K}})$ is a
Banach space. We define the map $\Lambda$ on $\mathcal{K}$ by
$$
\Lambda(v)=z,
$$
where $z$ is the mild solution to \eqref{sde1}.

\noindent We shall prove that $\Lambda$ is a contraction on
$\mathcal{K}$. For $i=1,\:2$, let $v_i$ in $\mathcal{K}$ and $z_i$
the solution to \eqref{sde1} corresponding to $v_i$, $i=1,\:2$. For
$n\geq 1$ we denote by $A_n$ the Yosida-approximation corresponding
to $A$. It is well known that for $n\geq 1$

$$
A_n=AJ_n\qquad \mbox{where}\quad J_n\defeq n(n-A)^{-1}.
$$

\noindent We now consider the approximating equation

\begin{equation}\label{appsde1}
 \left\{
\begin{array}{ll}
dz_n(t)=&\Big(A_nz_n(t)+F(z_n(t))\Big)dt+J_nB(v(t))dW_t,\quad
t\in [0,T],\\
z_n(0)\:=&x \in E,
\end{array}
\right.
\end{equation}

\noindent and let $z_i^n$ be the strong solution to \eqref{appsde1}
corresponding to $v_i$, $i=1,\:2$. (The $z_i^n$ are strong solutions
since $A_n$ is a bounded operator.)

 \noindent Hence by It\^{o}'s formula we have
\begin{equation}
\begin{split}
 \frac 1p\|z^n_1(t)-z^n_2(t)\|_0^p&=\int_0^t \|z^n_1(s)-z^n_2(s)\|_0^{p-2}\la A_n
(z^n_1(s)-z^n_2(s)),z^n_1(s)-z^n_2(s)\ra\:ds\\&+\int_0^t
\|z^n_1(s)-z^n_2(s)\|_0^{p-2} \la
F(z^n_1(s))-F(z^n_2(s)),z^n_1(s)-z^n_2(s)\ra\:ds\\&+\int_0^t
\|z^n_1(s)-z^n_2(s)\|_0^{p-2}
\|J_nB(v_1(s))-J_nB(v_2(s))\|^2_{\mathcal{L}_{HS}(U,H)}\:ds+M(t),
\end{split}
\end{equation}
\noindent where $$ M(t)\defeq \int_0^t \|z^n_1(s)-z^n_2(s)\|_0^{p-2}
\la z^n_1(s)-z^n_2(s),(J_nB(v_1(s))-J_nB(v_2(s)))\:dW(s)\ra.$$

\noindent By recalling the following inequality

\begin{equation*}
\la A n(n-A)^{-1}x,x\ra\leq \la A n(n-A)^{-1}x,
n(n-A)^{-1}x\ra,\quad x\in D(A),
\end{equation*}

\noindent using the definition of $A_n$, hypotheses ${\bf
(H_0)}$ and ${\bf (H_3)}$ it follows that
\begin{equation}\label{U}
\begin{split}
\frac 1p \|z^n_1(t)-z^n_2(t)\|_0^p&\leq -\omega\int_0^t
\|n(n-A)^{-1}(z_1^n(s)-z_2^n(s))\|_0^p\:ds +\kappa\int_0^t
\|z^n_1(s)-z^n_2(s)\|_0^p\:ds\\&+\int_0^t
\|z^n_1(s)-z^n_2(s)\|_0^{p-2}\|B(v_1(s))-B(v_2(s))\|^2_{\mathcal{L}_{HS}(U,H)}\:ds+M(t)\\
&\leq  \kappa\int_0^t \|z^n_1(s)-z^n_2(s)\|_0^p\:ds\\&+\int_0^t
\|z^n_1(s)-z^n_2(s)\|_0^{p-2}
\|B(v_1(s))-B(v_2(s))\|^2_{\mathcal{L}_{HS}(U,H)}\:ds+M(t)\\
&\leq  \kappa\int_0^t \|z^n_1(s)-z^n_2(s)\|_0^p\:ds\\&+\frac{p-2}{p}
\int_0^t \|z^n_1(s)-z^n_2(s)\|_0^p\:ds+ \frac{2}{p}L^{p}\int_0^t
\|v_1(s)-v_2(s)\|_0^p\:ds+M(t).
\end{split}
\end{equation}

\noindent This yields
\begin{equation*}
\begin{split}
\E\|z^n_1(t)-z^n_2(t)\|_0^p\leq \Big(p(\kappa+1)-2\Big)t\sup _{0\le
s\le t}\E \|z^n_1(s)-z^n_p(s)\|_0^p+ 2L^pt\sup _{0\le s\le t}\E
\|v_1(s)-v_2(s)\|_0^p.
\end{split}
\end{equation*}

\noindent Now by arguments similar to those of Proposition 7.17 and Theorem
7.18 in \cite{DPZ1} we have by letting $n\to +\infty$
\begin{equation*}
\begin{split}
\sup_{0\le t\le T} \E\|z_1(t)-z_2(t)\|_0^p &\leq
\Big(p(\kappa+1)-2\Big)T\sup _{0\le s\le t}\E
\|z_1(s)-z_2(s)\|_0^p\\&+ 2L^pT\sup _{0\le s\le t}\E
\|v_1(s)-v_2(s)\|_0^p.
\end{split}
\end{equation*}

\noindent Therefore, we have for $T$ small enough that
$$
\sup_{0\le t\le T} \E\|z_1(s)-z_2(s)\|_0^p\le \frac{1}{2} \sup_{0\le
t\le T} \E\|v_1(s)-v_2(s)\|_0^p.
$$

\noindent This shows that the mapping $\Lambda$ is a contraction on
$\mathcal{K}$ if $T$ is sufficiently small, and so it has a unique
fixed point $v$ in $\mathcal{K}$. The case of general $T>0$ can be
treated by considering the equation in intervals
$[0,\widetilde{T}]$, $[\widetilde{T},2\widetilde{T}]$, $\cdots$ for
small $\widetilde{T}$. The uniqueness follows by using the estimate
in \eqref{U} and taking expectation.
\end{proof}

\section{Invariant measures}
\noindent In this section we will prove existence of an invariant
measure $\mu$ for the process $\{u(t):\, t\geq 0\}$ given by
\eqref{sde0}. To this end we will use the Krylov-Bogoliubov Theorem.
So in particular we need to check tightness of the set of
probability measures $\left\{\mu_T\defeq
\frac{1}{T}\int_0^T\mu_{u(t,x)}\:dt,\;T\geq1\right\}$. Here
$\mu_{u(t,x)}$ denotes the distribution of $u(t,x)$,  $t\geq0$ with
$u(0)=x$. We remark that we will prove existence of an invariant
measure $\mu$ for \eqref{sde0} without condition on the size of the
Lipschitz constant $L$ of the diffusion term $B$. Therefore we need
an additional hypothesis on $F$.

\begin{enumerate}
 \item [${\bf (H_4)}$] There exists a continuous function $\rho:\:\R^+\rightarrow
 \R$, with $\lim\limits_{r\to+\infty}\frac{\rho(r^2)}{r^2}=-\infty$ such that
 $$
 \la F(u),u\ra\leq \rho(\|u\|^2_0),\quad u\in V_{\gamma}.
 $$
\end{enumerate}

\noindent  Note that hypothesis ${\bf (H_4)}$ implies that for all
$\lambda>0$ there exists $K_{\lambda}\geq 0$ such that

\begin{equation}\label{lambda}
\la F(v),v\ra\leq -\lambda\|v\|_0^2 +K_{\lambda}.
\end{equation}

%\noindent In fact, since $\lim\limits_{r\to+\infty}
%\frac{\rho(r^2)}{r^2}=-\infty$ we have for all $\lambda>0$ there
%exists $R_{\lambda}>0$ such that $\rho(\|v\|_0^2)\leq
%-\lambda\|v\|_0^2$, for any $v\in V_{\gamma}$ with $\|v\|_0^2\geq
%R_{\lambda}$. Hence using the continuity of $\rho$ we can write
%$$
%\rho(\|v\|_0^2)+\lambda\|v\|_0^2\leq K_\lambda,\quad\forall\: v\in
%V_{\gamma},
%$$

%for some  $K_{\lambda}\geq 0$. This yields the inequality
%\eqref{lambda}.

\noindent The following proposition shows tightness of the family of
measures $\left\{\mu_T,\;T\geq1\right\}$.
\begin{proposition}\label{tight}
Under hypotheses ${\bf(H_0)}$-${\bf(H_4)}$  the family of measures $
\left\{\mu_T,\: T\geq 1\right\}$ is tight.
\end{proposition}

\begin{proof}

\noindent Consider the solution $u(\cdot)$ of equation \eqref{sde0}.
If $(u(t))_{t\geq 0}$ is a strong solution (i.e, $u(t)\in D(A)$),
then by using It\^{o}'s formula and \eqref{lambda} we have for fixed
$t\geq 0$

\begin{equation}
\begin{split}
\E\|u(t)\|^2_0&= \E\|u(0)\|_0^2+2\E\int_0^t\la A(u(s)),u(s)\ra
ds\\&+2\E\int_0^t \la F(u(s)),u(s)\ra\:ds+\E\int_0^t\|B(u(s))\|_{\mathcal{L}_{HS}(U,H)}^2\:ds\\
&\leq \E\|u(0)\|_0^2+2
\E\int_0^t\Big(-c_\omega\|u(s)\|^{2}_{\gamma}-\lambda \|u(s)\|_0^2+
K_{\lambda}\Big)\:ds\\&\qquad+D\Big(t+\int_0^t\E\|u(s)\|_0^2\:ds\Big),
\end{split}
\end{equation}
where $D\defeq(L\vee\|B(0)\|_{\mathcal{L}_{HS}(U,H)})^2$ and
$c_{\omega}>0$, such that $c_{\omega}\|x\|_{\gamma}^2\leq
\|x\|_{\frac 12}^2$,  $x\in V_{\frac 12}$. In the case where
$u(\cdot)$ is a mild solution, starting in $u(0)=x\in E$, we shall consider the approximate equation
\[\left\{\begin{aligned}du_n(t)&=\left(A_n u_n(t)+F(u_n(t))\right)\,dt+J_n B(u_n(t))\,dW_t,\quad t\in [0,T],\\
u_n(0)&=x\in E,\end{aligned}\right.\]
compare with \eqref{appsde1}.

By It\^o's formula,
\begin{multline*}
\E\|u_n(t)\|_0^2=\|x\|_0^2+2\E\int_0^t\left\langle A_n(u_n(s)),u_n(s)\right\rangle\,ds\\
+2\E\int_0^t\left\langle F(u_n(s)),u_n(s)\right\rangle\,ds+\E\int_0^t\|J_n B(u_n(s))\|^2_{\mathcal{L}_{HS}(U,H)}\,ds\\
\le\E\|x\|_0^2+2\E\int_0^t\left(-\|(-A_n)^{\frac{1}{2}}u_n(s)\|_0^2-\lambda\|u_n(s)\|_0^2+K_\lambda\right)\,ds\\
+D\left(t+\E\int_0^t\|u_n(s)\|_0^2\,ds\right).
\end{multline*}
Pick $\lambda_\ast>0$ such that $\lambda_\ast>D/2$. Then we have
\[\E\|u_n(t)\|_0^2+(2\lambda_\ast-D)\E\int_0^t\|u_n(s)\|^2_0\,ds+2\E\int_0^t\|(-A_n)^{\frac{1}{2}}u_n(s)\|_0^2\,ds
\leq \|x\|^2_0+(D+2K_{\lambda_\ast})t.
\]
Now, by the results of Theorem \ref{existence}, we get that $u_n\to u$ in $C([0,T];L^2(\Omega;H))$.
Hence by Proposition \ref{prop:aux} in the Appendix and by Fatou's lemma,
\[\E\|u(t)\|_0^2+(2\lambda_\ast-D)\E\int_0^t\|u(s)\|^2_0\,ds+2\E\int_0^t\|(-A)^{\frac{1}{2}}u(s)\|_0^2\,ds
\leq \|x\|^2_0+(D+2K_{\lambda_\ast})t,
\]
which implies that $u(\cdot)\in L^2(\Omega\times [0,T];V_{\frac{1}{2}})$ and hence
\[\E\|u(t)\|_0^2+(2\lambda_\ast-D)\E\int_0^t\|u(s)\|^2_0\,ds+2c_\omega\E\int_0^t\|u(s)\|_\gamma^2\,ds
\le \|x\|^2_0+(D+2K_{\lambda_\ast})t.
\]

\noindent In particular, we have
$$
\E\frac{1}{t}\int_0^t\|u(s)\|^2_{\gamma}\:ds\leq
\frac{1}{2c_{\omega}}\Big(\E\|u(0)\|_0^2+2K_{\lambda_{\ast}}+D\Big)\quad\mbox{for
any $t\geq 1$}.
$$

\noindent We now take $\ve>0$, and put $R_{\ve}\defeq
\frac{1}{\sqrt{\varepsilon}}$. Then, for $T\geq 1$ we obtain
\begin{equation*}
\begin{split}
\mu_T(H\setminus\overline{B}(0,R_\ve)) & = \E\left(\frac
1T\int_0^T\ein_{\{\|u(s)\|_{\gamma}\geq R_\ve\}}\dd s\right)
    \leq\ve\E\left(\frac{1}{T}\int_0^T\|u(s)\|^2_{\gamma}\dd s\right)
\\
& \leq \ve
\frac{1}{2c_{\omega}}\Big(\E\|u(0)\|_0^2+2K_{\lambda_{\ast}}+D\Big).
\end{split}
\end{equation*}

\noindent Here, $\overline{B}(0,R_\ve)$ denotes the closed ball of
radius $R_\ve$ in $V_{\gamma}$. Since the embedding
$V_{\gamma}\hookrightarrow H$ is compact, the family of probability
measures $\{\mu_T \}_{T\ge 1}$ is tight on $H$. This completes the
proof.
\end{proof}

\noindent Now in order to conclude the existence of an invariant
measure for equation \eqref{sde0} we need to prove the Feller
property of $(u(t))_{t\geq 0}$.
\begin{proposition}
Assume hypotheses ${\bf(H_0)}$,  ${\bf(H_1)}$, ${\bf(H_2)}$ and
${\bf(H_3)}$. Let $(x_m)_{m\in\NN}$ be a sequence in $H$ such that
$x_m\xrightarrow[m\to +\infty]{\|\cdot\|_0} x$. Let $u^m$ (resp.
$u$) be the solutions to \eqref{sde0} with initial condition $x_m$
(resp. $x$). Then for any $t>0$,
\begin{equation}
\E\|u^m(t)-u(t)\|_0^2\rightarrow 0\quad \mbox{as}\:\:m\to+\infty.
\end{equation}

\noindent In particular, $(u(t))_{t\geq 0}$ is a Feller process.
\end{proposition}
\begin{proof} Assume that there is a strong solution $(u(t))_{t\geq 0}$, (i.e, $u(\cdot)\in D(A)$ and proceed by using Yosida-approximation for the general case. By using It\^{o}'s formula we obtain
\begin{equation}
\begin{split}
\E\|u^m(t)-u(t)\|_0^2&
\leq\E\|x-x_m\|_0^2+2\E\int_0^t\:\la
A(u^m(s))-A(u(s)),u^m(s)-u(s)\ra ds\\&+2\E\int_0^t\la
F(u^m(s))-F(u(s)),u^m(s)-u(s)\ra\:ds\\&\qquad +\E\int_0^t
\|B(u^m(s))-B(u(s))\|^2_{\mathcal{L}_{HS}(U,H)}\:ds\\
&\leq \E\|x-x_m\|_0^2+2(\kappa-\omega)\int_0^t\|u(s)-u^m(s)\|_0^2
ds+L\int_0^t \|u(s)-u^m(s)\|_0^2\:ds.
\end{split}
\end{equation}

\noindent Hence, by Gronwall's inequality,
\begin{equation}\label{invbound}
\E\|u^m(t)-u(t)\|_0^2\leq \|x_m-x\|_0^2\;\ee^{(2(\kappa-\omega)+L)t}.
\end{equation}

\noindent This implies in particular that for
$\psi:\:H\rightarrow\R$ bounded and continuous we have
$$
\lim\limits_{m\to+\infty}\E\psi(u^m(t))=\E\psi(u(t))\quad\mbox{for
any $t>0$},
$$
which yields the Feller property.
\end{proof}

\noindent Now, by the Krylov-Bogoliubov Theorem (see Section 3.1 in
\cite{DaZa:96}) we have the following result.
\begin{theorem}\label{inv}
Under hypotheses ${\bf(H_0)}-{\bf(H_4)}$ equation
\eqref{sde0} has an invariant measure.
\end{theorem}

\begin{rem}\label{inv2}
Assume hypotheses ${\bf(H_0)}-{\bf(H_4)}$.
Assume also that
\begin{equation}\label{eq:omega_large}\omega>\frac{L}{2}+\kappa.\end{equation}
 Then equation
\eqref{sde0} has a unique, ergodic, strongly mixing invariant measure.
\end{rem}
\begin{proof}
Taking \eqref{invbound} and \eqref{eq:omega_large} into account, the
claim follows by standard arguments.
See e.g.
\cite[proof of Proposition 2.2]{BDP}.
\end{proof}

\section{Applications}
Let $I = [0, L]\subset\R$ be a bounded interval and $A =
\frac{d^2}{dx^2}$ be the Laplacian with Dirichlet boundary
conditions. Clearly, $A$ is a negative definite self-adjoint
operator on $H = L^2 (I)$. The functions
$$
e_n (x) = \sqrt{\frac 2L} \sin\left(\frac{n\pi}{L} x\right) \, ,
n\ge 1 \, ,
$$
form an orthonormal set of eigenfunctions of $A$ with eigenvalues
$\lambda_n = -\left( \frac{\pi }{L}\right)^2 n^2$. For $\gamma>\frac
14$, we set $V_{\gamma}\defeq D((-\Delta)^{\gamma})$ and we define
$E\defeq C_0(I,\R)$ to be the Banach space of continuous real valued
functions on $I$ and vanishing at the boundary.

\noindent Let
\begin{equation}
\label{polynom} f(t) = a_{2n+1} t^{2n+1} + \ldots + a_1 t
\end{equation}
be a polynomial of odd degree with leading negative coefficient
$a_{2n+1} < 0$ and take $B$ a globally Lipschitz map from $H$ into
$\mathcal{L}_{HS}(H)$. We are interested in the stochastic partial
differential equation
\begin{equation}
\label{StochReactionEqu}\left\{\begin{array}{ll} du(t,x)=& \left(
\frac{d^2u}{dx^2}(t,x) + f(u(t,x))\right)\, dt + B(u(t,x))dW_t  \, ,
\quad (t,x) \in \R_+\times I\, ,\\
u(0,x)=&u_0(x),\quad u_0\in E.\end{array}\right.
\end{equation}
where $(W_t)_{t\geq0}$ is a cylindrical Wiener process on $L^2 (I)$.
For $u\in E$ define
$$
F(u) (x) = f(u(x)) \, , u\in E\:.
$$
\noindent Clearly $F$ maps $E$ into $E$ and is locally Lipschitz
continuous and bounded on bounded sets of $E$ and  by the Sobolev's
embedding theorem, the embedding  $V_{\gamma}\hookrightarrow E$ is
continuous for $\gamma>\frac 14$. Furthermore it is well known that
the part of the operator $A$ in $E$ generates a strongly continuous
semigroup on $E$. Hence hypothesis ${\bf (H_1)}$ is satisfied. By
using a characterization of the subdifferential of the norm in $E$
(see \cite[Example D.3]{DPZ1} it is not difficult to check
hypothesis ${\bf ( H_2)}$. Let us check hypothesis ${\bf (H_3)}$. We
can write
$$
F(u)=G_1(u)+G_2(u),
$$
where $G_1$ is dissipative (i.e., $\la G_1(u)-G_1(v),u-v\ra \leq
0,\:u,\:v\in E$) and $G_2$ Lipschitz continuous and bounded on $H$.
Indeed, let $\zeta_1$, $\zeta_2\in \R$ with $\zeta_1\leq \zeta_2$
such that $f(\zeta_1)>f(\zeta_2)$ and $f$ is decreasing on
$(-\infty, \zeta_1]\cup[\zeta_2,+\infty)$. Then by setting
$$
g_1(\zeta)=\left\{\begin{array}{ll}  f(\zeta),&\quad \zeta\in
(-\infty,
\zeta_1]\cup[\zeta_2,+\infty)\\
\ell(\zeta),&\quad \zeta\in [\zeta_1,\zeta_2], \end{array}\right.
$$
and $$ g_2(\zeta)=\left\{\begin{array}{ll}  0,&\quad \zeta\in
(-\infty,
\zeta_1]\cup[\zeta_2,+\infty)\\
f(\zeta)-\ell(\zeta),&\quad \zeta\in [\zeta_1,\zeta_2],
\end{array}\right.
$$
where
$\ell(\zeta)=f(\zeta_1)(\zeta_2-\zeta_1)^{-1}(\zeta_2-\zeta)+f(\zeta_2)(\zeta_2-\zeta_1)^{-1}(\zeta-\zeta_1)$
(the line which joins the points $(\zeta_1,f(\zeta_1)$ and
$(\zeta_2,f(\zeta_2)$), and defining
$$
G_1(u) (x) = g_1(u(x)),\qquad G_2(u) (x) = g_2(u(x)) \, , u\in E\, ,
$$

\noindent we see that $G_1$ and $G_2$ have the required properties.
Indeed, clearly $G_2$ is Lipschitz and bounded. For $G_1$, let $u$, $v$ in $E$
and set

$$
\Omega_u^1\defeq \{x\in I,\:\:u(x)\in [\zeta_1,\zeta_2]\},\quad
\Omega_v^1\defeq \{x\in I,\:\:v(x)\in [\zeta_1,\zeta_2]\},
$$

and
$$
\Omega_u^2\defeq \{x\in I,\:\:u(x)\in (-\infty,
\zeta_1)\cup(\zeta_2,+\infty)\},\quad \Omega_v^2\defeq \{x\in
I,\:\:v(x)\in (-\infty, \zeta_1)\cup(\zeta_2,+\infty)\}.
$$

\noindent Then
\begin{equation*}
\begin{split}
\la G_1(u)-G_1(v),u-v\ra&=\int_I (G_1(u(x))-G_1(v(x)))\cdot
(u(x)-v(x))\:dx\\
&=\int_{I\cap\Omega_u^1\cap\Omega_v^1} (G_1(u(x))-G_1(v(x)))\cdot
(u(x)-v(x))\:dx\\+& \int_{I\cap\Omega_u^1\cap\Omega_v^2}
(G_1(u(x))-G_1(v(x)))\cdot (u(x)-v(x))\:dx\\
&+\int_{I\cap\Omega_u^2\cap\Omega_v^1} (G_1(u(x))-G_1(v(x)))\cdot
(u(x)-v(x))\:dx\\+& \int_{I\cap\Omega_u^2\cap\Omega_v^2}
(G_1(u(x))-G_1(v(x)))\cdot (u(x)-v(x))\:dx\\ \le&
\int_{I\cap\Omega_u^1\cap\Omega_v^2}
(\ell(u(x))-f(v(x)))\cdot (u(x)-v(x))\:dx\\
&+\int_{I\cap\Omega_u^2\cap\Omega_v^1} (f(u(x))-\ell(v(x)))\cdot
(u(x)-v(x))\:dx\\+& \int_{I\cap\Omega_u^2\cap\Omega_v^2}
(f(u(x))-f(v(x)))\cdot (u(x)-v(x))\:dx\\
\end{split}
\end{equation*}

\noindent Clearly  $\int_{I\cap\Omega_u^2\cap\Omega_v^2}
(f(u(x))-f(v(x)))\cdot (u(x)-v(x))\:dx\leq 0$, since $f$ is
decreasing on $(-\infty, \zeta_1]\cup[\zeta_2,+\infty)$. On the
other hand for $x\in \Omega_u^1\cap\Omega_v^2$ it is not difficult to
see that

$$
(\ell(u(x))-f(v(x)))\cdot (u(x)-v(x))\leq 0.
$$

\noindent Similarly, in case $x\in \Omega_u^2\cap\Omega_v^1$ we have

$$(f(u(x))-\ell(v(x)))\cdot (u(x)-v(x))\leq 0.$$

\noindent This yields the required property for $G_1$ and therefore
hypothesis ${\bf (H_3)}$ is satisfied. Let us now prove hypothesis
${\bf (H_4)}$. To this end we write

\begin{equation*}
\begin{split}
\la F(u),u\ra=a_{2n+1}\int_0^1
u^{2n+2}(r)\:dr+\sum\limits_{k=1}^{2n}a_k\int_0^1 u^{k+1}(r)\:dr.
\end{split}
\end{equation*}

\noindent By using Young's inequality $ab\leq \frac\ve p a^p+\frac 1
{q \ve^{q-1}} b^q$, $p,\:q>1$, $pq=p+q$, $\ve>0$, we have
$$
\left|\int_0^1 u^{k+1}(r)\:dr\right|\leq \ve
\frac{k+1}{2n+2}\int_0^1u^{2n+2}(r)\:dr+ \ve^{-\frac{k+1}{2n-k+1}}.
$$
\medskip

\noindent Thus we can find some positive constant $C$ such that
$$
\la F(u),u\ra\leq \frac{a_{2n+1}}{2}\int_0^1u^{2n+2}(r)\:dr+ C.
$$

\noindent Since $a_{2n+1}<0$ we have
$$
\frac{a_{2n+1}}{2}\int_0^1u^{2n+2}(r)\:dr\leq \frac{a_{2n+1}}{2}
\|u\|^{2n+2}_{0}.
$$

\noindent Therefore, if we set $\rho(r)\defeq
\frac{a_{2n+1}}{2}r^{n+1}+C$, $r\in \R^+$ we have clearly
$\lim\limits_{r\to+\infty}\frac{\rho(r^2)}{r^2}=-\infty$ and
$$
\la F(u),u\ra \leq \rho(\|u\|^2_0),\quad u\in V_{\gamma} .
$$
\noindent This yields hypothesis ${\bf (H_4)}$.

\noindent By applying now Theorems \ref{existence} and \ref{inv} we
deduce that equation \eqref{StochReactionEqu} has a global solution
which belongs to $E$ and that \eqref{StochReactionEqu} has an
invariant measure.

\appendix
\section{$\Gamma$-convergence}

\begin{defi}
Let $q_n:H\to[0,+\infty]$, $n\in\mathbb{N}$, $q:H\to[0,+\infty]$ be closed, quadratic forms
(i.e. $q_n$, $q$ resp. have closed sublevel sets in $H$) with $q_n\not\equiv+\infty$, $n\in\mathbb{N}$,
$q\not\equiv+\infty$. We say that $\{q_n\}$ \emph{$\Gamma$-converges to} $q$ if the following holds true:\\
For $x_n\in H$, $n\in\mathbb{N}$, $x\in H$ such that $\|x_n-x\|_0\to 0$ as $n\to\infty$ it holds that
\begin{equation}\label{gamma1}\liminf_{n\to\infty}q_n(x_n)\ge q(x).\end{equation}
For each $y\in H$ there exist $y_n\in H$, $n\in\mathbb{N}$, with $\|y_n-y\|_0\to 0$ as $n\to\infty$ and
\begin{equation}\label{gamma2}\limsup_{n\to\infty}q_n(y_n)\le q(y).\end{equation}
\end{defi}

\begin{proposition}\label{prop:aux}
Let $A$ be as in the main part and let $A_n:=nA(n-A)^{-1}$ be its Yosida approximation. Let
\[\Phi_n(u):=\|(-A_n)^{\frac{1}{2}}u\|_0^2,\quad u\in H,\]
furthermore, let
\[\Phi(u):=\|(-A)^{\frac{1}{2}}u\|_0^2,\quad u\in D((-A)^\frac{1}{2}).\]
Extend $\Phi$ to $H$ by $\Phi(u):=+\infty$ whenever $u\in H\setminus D((-A)^{\frac{1}{2}})$.

Then $\{\Phi_n\}$ $\Gamma$-converges to $\Phi$.
\end{proposition}
\begin{proof}
First observe that $\Phi$ is a closed quadratic form on $H$ associated to the positive self-adjoint operator $-A$, see \cite[Chapter 12]{DaMa}. By \cite[Proposition 12.23]{DaMa}, $\Phi_n$ equals the
so-called
Moreau-Yosida approximation
\[\inf_{y\in H}\left[\Phi(y)+n\|y-x\|^2_0\right]\]
of $\Phi$. By \cite[Theorem 9.13, Corollary 9.14]{DaMa}, we see that $\Phi_n\uparrow\Phi$
pointwise as $n\to\infty$. The claim follows now by \cite[Remark 5.5]{DaMa}.
\end{proof}

\end{document}